\documentclass[12pt]{amsart}
\usepackage{amsmath,amsfonts,amssymb,a4wide,url}
\newcommand{\eps}{\varepsilon}

\theoremstyle{plain}
\newtheorem{lemma}{Lemma}
\newtheorem{theorem}{Theorem}
\newtheorem{corollary}[theorem]{Corollary}

\theoremstyle{definition}
\newtheorem{defn}{Definition}
\newtheorem{remark}{Remark}

\newtheorem{example}{Example}
\newtheorem*{ackno}{Acknowledgment}
\newcommand{\Z}{\mathbb{Z}}

\newcommand{\C}{\mathbb{C}}
\newcommand{\N}{\mathbb{N}}
\newcommand{\bfell}{\boldsymbol{\ell}}
\newcommand{\BigOh}{\mathcal{O}}

\newcommand{\fl}[1]{\lfloor #1\rfloor}
\DeclareMathOperator*{\Res}{\mathrm{Res}}

\numberwithin{equation}{section}

\begin{document}
\title[The Zeta Function\ldots]
{The Zeta Function of the Laplacian on Certain Fractals}

\author[G. Derfel]{Gregory Derfel}
\address[G. D.]{Department of Mathematics and Computer Science,
Ben Gurion University of the Negev,
Beer Sheva 84105,
Israel}
\email{derfel@math.bgu.ac.il} 

\author[P. J. Grabner]{Peter J. Grabner\dag{}}
\thanks{\dag{} This author is supported by the START project Y96-MAT of the
  Austrian Science Fund}
\address[P. G.]{Institut f\"ur Mathematik A,
Technische Universit\"at Graz,
Steyrergasse 30,
8010 Graz,
Austria}
\email{peter.grabner@tugraz.at}

\author[F. Vogl]{Fritz Vogl}
\address[F. V.]{Institut f\"ur Analysis und Scientific Computing,
Technische Universit\"at Wien,
Wiedner Hauptstra\ss{}e 8--10,
1040 Wien,
Austria}

\email{fvogl@osiris.tuwien.ac.at}
\keywords{Dirichlet series, Laplace operator, fractals, spectral decimation,
  complex dynamics}
\subjclass[2000]{Primary: 30B50; Secondary: 11M41 37F10}
\date{\today}
\begin{abstract}
  We prove that the zeta-function $\zeta_\Delta$ of the Laplacian $\Delta$ on a
  self-similar fractals with spectral decimation admits a meromorphic
  continuation to the whole complex plane. We characterise the poles, compute
  their residues, and give expressions for some special values of the
  zeta-function. Furthermore, we discuss the presence of oscillations in the
  eigenvalue counting function.
\end{abstract}

\maketitle

\section{Introduction}\label{sec:introduction}
Diffusion on fractals has been studied extensively as a generalisation of usual
Brownian motion on manifolds. After its introduction in the physics literature
(cf.~\cite{Rammal_Toulouse1983:random_walks_fractal}) M.~Barlow and E.~Perkins
\cite{Barlow_Perkins1988:brownian_motion_sierpinski} gave a very detailed study
of the properties of the diffusion on the Sierpi\'nski gasket. Later
T.~Lindstr\o{}m \cite{Lindstroem1990:brownian_motion_nested} generalised these
results to nested fractals. Especially, he derived results on the distribution
of the eigenvalues of the Laplacian associated to the diffusion. Here the
Laplacian is seen as the infinitesimal generator of Brownian motion.
Alternatively, the Laplacian can be obtained as the limit of properly rescaled
difference operators on graphs approximating the fractal structure
(cf.~\cite{Fukushima_Shima1992:spectral_analysis_sierpinski,
  Kigami1993:harmonic_calculus_p,Kigami2001:analysis_fractals}). Furthermore,
the Laplace operator can be obtained via Dirichlet forms
(cf.~\cite{Barlow1998:diffusions_on_fractals}).

There exists a vast literature on properties of the spectrum of the Laplacian.
For instance, we refer to
\cite{Barlow_Kigami1997:localized_eigenfunctions_laplacian,
  Kigami1998:distributions_localized_eigenvalues, Kigami2001:analysis_fractals,
  Kigami_Lapidus1993:weyl_problem_spectral, Kroen2002:green_functions_self,
  Lapidus1994:analysis_fractals_laplacians,
  Malozemov_Teplyaev1995:pure_point_spectrum, Sabot2000:pure_point_spectrum,
  Shima1993:eigenvalue_problem_laplacian,
  Shima1996:eigenvalue_problems_laplacians,
  Strichartz1999:some_properties_laplacians}. Especially, it has been proved
that the Laplacian on an infinite post-critically-finite (p.~c.~f.)
self-similar fractal has pure point spectrum
(cf.~\cite{Sabot2000:pure_point_spectrum}).

In the case of a compact self-similar fractal $G$ equipped with the Hausdorff
measure $\mathcal{H}$ it is natural to ask for the behaviour of the counting
function of the eigenvalues under Dirichlet or Neumann boundary conditions
\begin{equation}\label{eq:count}
L(x)=\sum_{\substack{-\Delta u=\mu u\\\mu<x}}1,
\end{equation}
the trace of the heat kernel $p_t(x,y)$
\begin{equation}\label{eq:heat}
P(t)=\sum_{-\Delta u=\mu u}e^{-\mu t}=\int_Gp_t(x,x)\,d\mathcal{H}(x),
\end{equation}
or the corresponding Dirichlet generating function
\begin{equation}\label{eq:zeta}
\zeta_\Delta(s)=\sum_{-\Delta u=\mu u}\frac1{\mu^s},
\end{equation}
the ``zeta-function'' of the Laplace operator (the eigenvalues are counted
\emph{with} multiplicities in all three sums). This zeta-function will be the
subject of the present paper.

In the case of the Laplacian on a Riemannian manifold, the meromorphic
continuation of $\zeta_\Delta$ can be derived directly from the asymptotic
expansion of the trace of the heat-kernel
(cf.~\cite{Minakshisundaram_Pleijel1949:some_properties_eigenfunctions,
  Rosenberg1997:laplacian_riemannian_manifold}). In this case it is known that
$\zeta_\Delta$ admits a meromorphic continuation to the whole complex plane with
at most simple poles at $\frac d2,\frac d2-1,\ldots,\frac
d2-\lfloor\frac{d-1}2\rfloor$, where $d$ denotes the dimension of the manifold.
Furthermore, the residues at these poles have a geometric meaning.

In the case of a fractal the situation is more delicate. The current knowledge
on the behaviour of the heat-kernel is not sufficient to provide enough
information on its trace to get the meromorphic continuation of $\zeta_\Delta$.
In the present paper we will present a method for the continuation of
$\zeta_\Delta$ which is based on an asymptotic study of the solution of
Poincar\'e's functional equation. In Section~\ref{sec:asympt} we prove several
results relating the behaviour of the solution to properties of the according
Julia set. This approach for the analytic continuation of $\zeta_\Delta$ is
restricted to a subclass of self-similar fractals which exhibit spectral
decimation (cf.~\cite{Rammal_Toulouse1983:random_walks_fractal}) with a
polynomial map (see Section~\ref{sec:eigenv-lapl-spectr}). On the other hand it
allows to attach a zeta function to every polynomial which has a certain
behaviour under iteration. Such zeta-functions have been introduced in
\cite{Teplyaev2005:spectral_zeta_functions}; we give a totally different
approach for their meromorphic continuation and are also able to derive
expressions for some special values.

In Section~\ref{sec:analyt-cont-zeta} we use the properties of the zeta
function to confirm a conjecture
(cf.\cite[p.~105]{Kigami_Lapidus1993:weyl_problem_spectral}) about the presence
of oscillations in the eigenvalue counting function of the Laplacian for the
class of fractals under consideration.
\section{Eigenvalues of the Laplacian and spectral decimation}
\label{sec:eigenv-lapl-spectr}
As it was observed in the case of the Sierpi\'nski gasket and its higher
dimensional analogues in
\cite{Fukushima_Shima1992:spectral_analysis_sierpinski,
  Shima1993:eigenvalue_problem_laplacian,
  Shima1996:eigenvalue_problems_laplacians} the eigenvalues of the Laplace
operator on certain finitely ramified self-similar fractals exhibit the
phenomenon of spectral decimation
(cf.~\cite{Fukushima_Shima1992:spectral_analysis_sierpinski} for the case of
higher-dimensional Sierpi\'nski gaskets and
\cite{Malozemov_Teplyaev2003:self_similarity_operators,
  Strichartz2003:fractafolds_spectra, Teplyaev2005:spectral_zeta_functions} for
a more general treatment).
\begin{defn}[Spectral decimation]\label{def1}
  The Laplace operator on a p.~c.~f.\ self-similar fractal $G$ admits
  \emph{spectral decimation}, if there exists a rational function $R$, a finite
  set $A$ and a constant $\lambda>1$ such that all eigenvalues of $\Delta$ can
  be written in the form 
\begin{equation}\label{eq:decimation}
\lambda^m\lim_{n\to\infty}\lambda^nR^{(-n)}(\{w\}),\quad w\in A, m\in\N
\end{equation}
where the preimages of $w$ under $n$-fold iteration of $R$ have to be chosen
such that the limit exists. Furthermore, the multiplicities of the
eigenvalues depend only on $w$ and $m$, and the generating functions of the
multiplicities are rational.
\end{defn}
The fact that all eigenvalues
of $\Delta$ are negative real implies that the Julia set of $R$ has to be
contained in the negative real axis. We will exploit this fact later.

In many cases such as the higher dimensional Sierpi\'nski gaskets, the rational
function $R$ is a polynomial. The method for meromorphic continuation of
$\zeta_\Delta$ given in Section~\ref{sec:weierstr-prod-dirich} makes use of
this assumption.  In a recent paper \cite{Teplyaev2005:spectral_zeta_functions}
A.~Teplyaev showed under the same assumption that the $\zeta$-function of the
Laplacian admits a meromorphic continuation to $\Re s>-\eps$ for some $\eps>0$
depending on properties of the Julia set of the polynomial given by spectral
decimation. His method uses ideas similar to those used in
\cite{Grabner1997:functional_iterations_stopping} for the meromorphic
continuation of a Dirichlet series attached to a polynomial. Complementary to
the ideas used here, Teplyaev's method carries over to rational functions $R$.

\section{Asymptotic behaviour of solutions of a functional equation}
\label{sec:asympt}
In this section we collect several results on the growth of the solution of the
Poincar\'e functional equation related to polynomials. We note here that
Valiron \cite[Chapter~VII]{Valiron1954:fonctions_analytiques} gave rather
general results on the order of the solution of such functional equations.
Throughout the paper we will use the notation
\begin{equation*}
\log_\lambda z=\frac{\log z}{\log\lambda}
\end{equation*}
for the logarithm to base $\lambda$. We first give a coarse growth estimate,
which is contained in \cite{Valiron1954:fonctions_analytiques} and will be
refined in Theorem~\ref{thm1}.
\begin{lemma}\label{lem0}
  Let $p(x)=a_dx^d+\cdots+a_1x$ be a polynomial with real coefficients,
  $d\geq2$, which satisfies $p(0)=0$ and $p'(0)=\lambda>1$. Then the solution
  $\Phi$ of the functional equation
\begin{equation}\label{eq:schr}
\Phi(\lambda z)=p(\Phi(z)),\quad \Phi(0)=0,\quad\Phi'(0)=1
\end{equation}
is an entire function of order $\leq\rho=\log_\lambda d$.
\end{lemma}
\begin{proof}  By the general theory of iteration of polynomial functions,
  there exists a unique entire solution $\Phi$ of the equation \eqref{eq:schr}
  (cf.~\cite{Beardon1991:iteration_rational_functions,
    Kuczma1963:schroeder_equation}).  Furthermore, there exists a constant
  $C>1$, such that $|p(z)|\leq C\max(1,|z|^d)$. Therefore,
\begin{equation*}
|\Phi(z)|\leq C^{d^{n-1}+\ldots+d+1}\max(1,|\Phi(\lambda^{-n}z)|^{d^n}).
\end{equation*}
Let $r>0$ such that $|\Phi(z)|\leq1$ for $|z|\leq r$ and take 
$n=[\log_\lambda(|z|/r)]+1$ to obtain
\begin{equation*}
|\Phi(z)|\leq C^{dr^{-\rho}|z|^\rho}.
\end{equation*}
\end{proof}
The following theorem gives a more precise description of the growth order of
$\Phi$ under assumptions on the behaviour of $p$ under iteration.
\begin{theorem}\label{thm1}
  Let $p(x)=a_dx^d+\cdots+a_1x$ be a polynomial with real coefficients,
  $d\geq2$, which satisfies $p(0)=0$ and $p'(0)=\lambda>1$.  Furthermore,
  suppose that $\mathcal{F}_\infty$, the Fatou component of $\infty$  of $p$,
  contains an angular region of the form
\begin{equation*}
W_\beta=\left\{z\in\mathbb{C}\setminus\{0\}\mid |\arg z|<\beta\right\}
\end{equation*}
for some $\beta>0$.
Then the solution $\Phi$ of the functional equation \eqref{eq:schr}
is an entire function. For any $\eps>0$ and any $M>0$ the asymptotic relation
\begin{equation}\label{eq:Phiasymp}
\Phi(z)=\frac1{\sqrt[d-1]{a_d}}\exp\left(z^\rho F(\log_\lambda z)+
o\left(|z|^{-M}\right)\right)
\end{equation}
holds uniformly for $|\arg z|\leq\beta-\eps$, where $\rho=\log_\lambda d$ and
$F$ is a periodic holomorphic function of period $1$ on the strip
$\{w\in\mathbb{C}\mid|\Im w|<\frac\beta{\log\lambda}\}$. The real part of
$z^\rho F(\log_\lambda z)$ is always positive; $F$ takes real values on the
real axis.
\end{theorem}
\begin{proof}
  After Lemma~\ref{lem0} it remains to prove only the assertion about the
  asymptotic behaviour of $\Phi$.
  
  Since $p$ is a polynomial of degree $d\geq2$, there exists a real number
  $R>0$, such that $|z|\geq R\Rightarrow|p(z)|\geq 2|z|$, which implies that
  $(p^{(n)}(z))_n$ tends to infinity for $|z|>R$. Fix $\eps>0$. Then by local
  conformity of $\Phi$ around $0$ there exists $r>0$ such that
  \begin{equation*}
    \Phi\left(\{z\in\mathbb{C}\setminus\{0\}\mid |z|\leq r,
 |\arg z|\leq\beta-\eps\} \right)
\subset W_\beta.
  \end{equation*}
By definition of the Fatou set the sequence $(p^{(n)}(z))_n$ is equicontinuous
on $W_\beta$ and therefore contains a subsequence $(p^{(n_k)}(z))_k$, which
converges to $\infty$ uniformly on 
\begin{equation*}
M_\eps=
\Phi\left(\{z\in\mathbb{C}\setminus\{0\}\mid \frac r\lambda\leq|z|\leq r,
 |\arg z|\leq\beta-\eps\} \right)
\end{equation*}
by the fact that $\infty$ is an attracting fixed point of $p$ and our
assumptions on the Fatou component of $\infty$. Thus there exists an $N$ such
that $|p^{(N)}(z)|>R$ for all $z\in M_\eps$, which implies that
$(p^{(n)}(z))_n$ converges to $\infty$ uniformly on $M_\eps$. Inserting this
fact into \eqref{eq:schr} this implies that
\begin{equation}\label{eq:Phiinfty}
\lim_{z\to\infty}|\Phi(z)|=\infty\text{ uniformly for }|\arg z|\leq\beta-\eps.
\end{equation}
Furthermore, for any $K>1$ there exists a $C>0$ such that
\begin{equation*}
\left|\Phi(z)\right|\geq C|z|^K\text{ for }|\arg z|\leq\beta-\eps.
\end{equation*}

Now note that our assumption on the Fatou set implies that there are no zeros
of $\Phi$ in $W_\beta$. Therefore,
\begin{equation*}
\Psi(z)=\log\Phi(z)
\end{equation*}
is analytic in $W_\beta$ and satisfies the functional equation
\begin{equation}\label{eq:Psi}
\Psi(\lambda z)=d\Psi(z)+\log a_d+r(z)\text{ with }
r(z)=\log\frac{p(\Phi(z))}{a_d\Phi(z)^d}.
\end{equation}
The function $r(z)$ is analytic in
\begin{equation*}
\left\{z\in\mathbb{C}\mid |z|>r_0, |\arg z|<\beta\right\}
\end{equation*}
for some $r_0>0$.
Setting $\Psi(z)=z^\rho\phi(z)-\frac{\log a_d}{d-1}$ with $\rho=\log_\lambda d$
yields
\begin{equation}\label{eq:phifunc}
\phi(\lambda z)=\phi(z)+\frac{r(z)}{dz^\rho}.
\end{equation}
By our previous knowledge on $\Phi(z)$ for $z\to\infty$ in $|\arg z|<\beta$ it
follows that $r(z)$ tends to $0$ faster than any negative power of $|z|$
uniformly in any angular region $|\arg z|\leq\beta-\eps$. Iterating
\eqref{eq:phifunc} yields
\begin{equation}\label{eq:phin}
\phi(\lambda^nz)=\phi(z)+
\frac1{z^\rho}\sum_{k=0}^{n-1}\frac{r(\lambda^kz)}{d^{k+1}}.
\end{equation}
We set
\begin{equation*}
F\left(\log_\lambda z\right)=\lim_{n\to\infty}\phi(\lambda^nz),
\end{equation*}
which exists by our knowledge on $r(z)$ and \eqref{eq:phin}. Clearly, $F(w)$ is
periodic with period $1$, holomorphic in $|\Im w|<\frac\beta{\log\lambda}$ and
\begin{equation*}
\phi(z)=F\left(\log_\lambda z\right)-\frac1{z^\rho}
\sum_{k=0}^\infty\frac{r(\lambda^kz)}{d^{k+1}}=
F\left(\log_\lambda z\right)+o(|z|^{-M})
\end{equation*}
for any positive $M$.

Putting everything together, we obtain
\begin{equation*}
\Psi(z)=z^\rho F\left(\log_\lambda z\right)-\frac{\log a_d}{d-1}+
o(|z|^{-M}),
\end{equation*}
which implies \eqref{eq:Phiasymp}. The statement on the real part of
$z^\rho F(\log_\lambda z)$ follows from $|\Phi(z)|\to\infty$ for $|z|\to\infty$
in $W_\beta$.
\end{proof}

\begin{corollary}\label{cor1}
The limit
\begin{equation*}
\lim_{n\to\infty}\frac{\log\Phi(\lambda^nz)}{d^nz^\rho}
\end{equation*}
exists for all $z\in\Phi^{(-1)}(\mathcal{F}_\infty)$. This limit gives the
analytic continuation of $F(\log_\lambda z)$ to
$\Phi^{(-1)}(\mathcal{F}_\infty)$.
\end{corollary}
\begin{proof}
The proof of the asymptotic formula for $\Psi(z)$ used only the fact  that
$\Phi(\lambda^nz)\to\infty$ for $n\to\infty$. This relation holds for any point
$z\in\Phi^{(-1)}(\mathcal{F}_\infty)$.
\end{proof}
\begin{remark}\label{rem1}
  Constancy of the periodic function $F$ occurring in the statement of
  Theorem~\ref{thm1} has been the subject of a series of papers in the context
  of branching processes \cite{Biggins_Bingham1991:near_constancy_phenomena,
    Dubuc1982:etude_theorique_numerique}, where it turns out that the constancy
  of $F$ is implied by the existence of a continuous time extension of the
  according branching process
  (cf.~\cite{Karlin_Mcgregor1968:embeddability_discrete_time}). Furthermore, it
  is conjectured that also the opposite implication is true. Usually, on the
  real axis the function $F$ exhibits very small oscillations around a mean
  value, which can only be observed using high precision numerical
  computations. Therefore, theoretical conditions for the presence of such
  oscillations are of special interest.
\end{remark}
For the purposes of this paper we give the following condition for
non-constancy of the periodic function $F$.
\begin{theorem}\label{thm2}
Let $\Phi$ be the solution of the functional equation \eqref{eq:schr}. Assume
that there exists an angle $\gamma$ such that $\Phi(\{re^{i\gamma}\mid r>0\})$
intersects the Fatou-component $\mathcal{F}_\infty$ as well as the Julia-set
$\mathcal{J}_p$. Then the periodic function $F$ in \eqref{eq:Phiasymp} is not
constant.
\end{theorem}
\begin{proof}
We note first that $\mathcal{J}_p$ is a compact subset of $\mathbb{C}$, since
$p$ is a polynomial
(cf.~\cite{Beardon1991:iteration_rational_functions}). Assume that
$F(\log_\lambda z)=C$. Then by Corollary~\ref{cor1} we have
\begin{equation*}
\lim_{n\to\infty}\frac{\log\Phi(\lambda^nz)}{d^nz^\rho}=C
\end{equation*}
for $z\in\Phi^{(-1)}(\mathcal{F}_\infty)$. For any $r>0$ with
$\Phi(re^{i\gamma})\in\mathcal{J}_p$ we have
$\Phi(r\lambda^ne^{i\gamma})\in\mathcal{J}_p$ for all $n$. Thus by our
assumptions, for any $k\in\mathbb{N}$, any $\eps>0$, and any $M>0$ there exists
$r>M$, such that $\Phi(re^{i\gamma})\in\mathcal{F}_\infty$ and
\begin{equation*}
\inf\{|\Phi(re^{i\gamma})-z|\mid z\in\mathcal{J}_p\}<\eps.
\end{equation*}
Then we have
\begin{align*}
&\lim_{n\to\infty}\frac{\log\Phi(\lambda^{n+k}re^{i\gamma})}
{r^\rho d^{n+k}e^{i\rho\gamma}}=C\\
&\lim_{n\to\infty}\frac{\log p^{(k)}(\Phi(\lambda^nre^{i\gamma}))}
{d^n(p^{(k)}(\Phi(re^{i\gamma})))^\rho}=C.
\end{align*}
On the other hand
\begin{equation*}
\lim_{n\to\infty}\frac{\log\Phi(\lambda^{n+k}re^{i\gamma})}
{r^\rho d^{n+k}e^{i\rho\gamma}}=
\lim_{n\to\infty}\frac{\log p^{(k)}(\Phi(\lambda^nre^{i\gamma}))}
{d^n(p^{(k)}(\Phi(re^{i\gamma})))^\rho}
\frac{(p^{(k)}(\Phi(re^{i\gamma})))^\rho}{d^k r^\rho e^{i\gamma\rho}}=C.
\end{equation*}
Since $p^{(k)}(\Phi(re^{i\gamma}))$ is contained in a bounded set and $r$ can
be chosen arbitrarily large. This gives a contradiction.
\end{proof}
\begin{corollary}\label{cor2}
  Assume that $\mathcal{J}_p\subset\mathbb{R}^-$ and that $\Phi$ is unbounded
  on $\mathbb{R}^-$. Then the periodic function $F$ in Theorem~\ref{thm1} is
  not constant.
\end{corollary}
\begin{proof}
  $\Phi(\mathbb{R}^-)$ intersects the Julia set $\mathcal{J}_p$, since $\Phi$
  attains negative values close to $0$, and $0\in\mathcal{J}_p$.
  $\Phi(\mathbb{R}^-)$ intersects $\mathcal{F}_\infty$, since $\Phi$ takes
  arbitrarily large values. Then Theorem~\ref{thm2} gives the assertion with
  $\gamma=\pi$.
\end{proof}
\begin{remark}\label{rem3}
  The example of $p(x)=x(4+x)$, where $\Phi(z)=4\sinh^2\frac12\sqrt{z}$, shows
  that the assumption on the unboundedness of $\Phi$ on $\mathbb{R}^-$ cannot
  be omitted.
\end{remark}
\begin{corollary}\label{cor3}
  If $\rho<\frac12$, then the periodic function $F$ in Theorem~\ref{thm1} is
  not constant.
\end{corollary}
\begin{proof}
  By \cite[Theorem~3.1.5]{Boas1954:entire_functions} a function of order
  $<\frac12$ is unbounded on any ray, and so is $\Phi$. Apply
  Corollary~\ref{cor2} to a ray, whose $\Phi$-image intersects the Julia set of
  $p$ to obtain the assertion.
\end{proof}
\begin{corollary}\label{cor4}
  If the periodic function $F$ in Theorem~\ref{thm1} is not constant, it has
  infinitely many non-zero Fourier-coefficients.
\end{corollary}
\begin{proof}
  Since $p(0)=0$ and $0$ lies on the boundary of $\mathcal{F}_\infty$ by the
  hypotheses of Theorem~\ref{thm1}, infinitely many zeros of $\Phi$ lie on the
  boundary of $\Phi^{(-1)}(\mathcal{F}_\infty)$. From this we get that
  $F(\log_\lambda z)$ is unbounded when approaching these boundary points of
  $\Phi^{(-1)}(\mathcal{F}_\infty)$, which implies that $F$ has to have
  infinitely many non-zero Fourier-coefficients.
\end{proof}
\section{Spectral decimation and level sets of $\Phi$}
\label{sec:spectr-decim-level}
By \eqref{eq:decimation} all eigenvalues of the Laplacian $\Delta$ can be
computed by iterating $p^{(-1)}$ and rescaling. In this section we want to
analyse this more precisely. The inverse map of the polynomial function $p$ has
$d=\deg p$ branches, which we denote by $q_1,\ldots,q_d$. We denote by $q_1$
the unique branch given by $q_1(0)=0$, which also satisfies
$q'_1(0)=\frac1\lambda$. Given $w\in\mathbb{C}$ the elements of
$p^{(-n)}(\{w\})$ are in correspondence to words of length $n$ over the
alphabet $\{1,\ldots,d\}$: associate to
$\ell=(\ell_1,\ldots,\ell_n)\in\{1,\ldots,d\}^n$ the value
$z_\ell=q_{\ell_n}\circ q_{n-1}\circ\cdots\circ q_{\ell_1}(w)$. Then
$p^{(n)}(z_\ell)=w$ and all elements of $p^{(-n)}(\{w\})$ are of this form.

Let $\bfell=(\ell_1,\ell_2,\ldots)\in\{1,\ldots,d\}^{\N}$. Then the limit
\begin{equation*}
\lim_{n\to\infty}\lambda^n q_{\ell_n}\circ 
q_{n-1}\circ\cdots\circ q_{\ell_1}(w)
\end{equation*}
exists if and only if there exists $N\in\N$, such that $\ell_n=1$ for $n\geq
N$. This follows from the fact that $q_1$ is the only branch of $p^{(-1)}$
which maps $0$ to itself.

Take $z_n\in\lambda^n p^{(-n)}(\{w\})$. Then define $u_n\in\C$ by
$\lambda^{-n}z_n=\Phi(\lambda^{-n}u_n)$. Then
$z_n=\lambda^n\Phi(\lambda^{-n}u_n)$ and assuming that the limit
$\lim_{n\to\infty}z_n$ exists as discussed above, we have 
\begin{equation*}
z=\lim_{n\to\infty}z_n=\lim_{n\to\infty}\lambda^n\Phi(\lambda^{-n}u_n)=
\lim_{n\to\infty}u_n.
\end{equation*}
On the other hand $w=p^{(n)}(\Phi(\lambda^{-n}u_n))=\Phi(u_n)$ and therefore
$\Phi(z)=w$. Thus the solutions of $\Phi(z)=w$ are exactly the possible limits
of sequences $(z_n)$ as above. Summing up, we have
\begin{lemma}\label{lem5}
Let $\bfell=(\ell_1,\ell_2,\ldots)\in\{1,\ldots,d\}^{\N}$ with
$\ell_N=\ell_{N+1}=\ldots=1$ for some $N\in\N$. Then the limit
\begin{equation*}
z_{\bfell}=\lim_{n\to\infty}\lambda^n
q_{\ell_n}\circ q_{n-1}\circ\cdots\circ q_{\ell_1}(w)
\end{equation*}
exists and satisfies $\Phi(z_{\bfell})=w$. All solutions of the equation
$\Phi(z)=w$ can be obtained as such a limit.
\end{lemma}
\section{Weierstrass products and Dirichlet series}
\label{sec:weierstr-prod-dirich}
Since the solution $\Phi$ of the functional equation \eqref{eq:schr} is of
order $\rho$, we can represent $\Phi_w(z)=1-\frac1w\Phi(z)$ by an Hadamard
product (for $w\neq0$)
\begin{multline}\label{eq:hadamard}
\Phi_w(z)=\exp(c_1(w)z+\cdots+c_k(w)z^k)\\
\times\prod_{\ell=1}^\infty\left(1+\frac{z}{\mu_\ell(w)}\right)
\exp\left(-\frac{z}{\mu_\ell(w)}+\frac{z^2}{2\mu_\ell(w)^2}+\cdots
+(-1)^k\frac{z^k}{k\mu_\ell(w)^k}\right),
\end{multline}
where $k=\fl{\rho}$ denotes the integer part of $\rho$
(cf.~\cite{Boas1954:entire_functions}). It is clear from the validity of the
functional equation \eqref{eq:schr} that $\Phi_w$ has infinitely many zeroes
even in the case of integer $\rho$.

It follows from Theorem~\ref{thm1} that $\log\Phi_w(z)$ has a convergent
Taylor series around $0$:
\begin{equation}\label{eq:logPhi}
\log\Phi_w(z)=\sum_{\ell=1}^\infty b_\ell(w) z^\ell.
\end{equation}
Furthermore, assuming that the Julia set of $p$ is a subset of the negative
reals and $w<0$ as it is the case for spectral decimation, $\log\Phi_w(z)$ is
holomorphic in $\{z\in\mathbb{C}\mid|\arg z|<\beta\}$. Computing the Taylor
series of $\log\Phi_w(z)$ around $z=0$ and comparing coefficients yields
$c_\ell(w)=b_\ell(w)$ for $\ell=1,\ldots,k$.

Let
\begin{equation*}
q_m(z,w)=\sum_{\ell=1}^m b_\ell(w) z^\ell
\quad\text{for }m\in\mathbb{N}
\end{equation*}
and consider the Mellin transform of $\log\Phi_w(x)-q_k(x,w)$
\begin{equation}\label{eq:mellin}
M_w(s)=\int\limits_0^\infty\left(\log\Phi_w(x)-q_k(x,w)\right)x^{s-1}\,dx.
\end{equation}
The function $M_w(s)$ is holomorphic in the strip $-k-1<\Re s<-\rho$ by
general properties of the Mellin transform
(cf.~\cite{Doetsch1971:handbuch_der_laplace}). In order to obtain a meromorphic
continuation of $M_w(s)$ we split the integral in \eqref{eq:mellin} into two
parts:
\begin{align}\label{eq:mellin1}
M_w^{(1)}(s)&=\int\limits_0^1\left(\log\Phi_w(x)-q_k(x,w)\right)x^{s-1}\,dx\\
M_w^{(2)}(s)&=\int\limits_1^\infty
\left(\log\Phi_w(x)-q_k(x,w)\right)x^{s-1}\,dx.
\label{eq:mellin2}
\end{align}

It is clear that $M_w^{(1)}(s)$ is holomorphic in $\Re s>-k-1$.  From the
expansion \eqref{eq:logPhi} of $\log\Phi_w(x)$ we immediately get a meromorphic
continuation of $M_w^{(1)}(s)$ into $\Re s>-r$ for $r>k+1$, $r\in\mathbb{N}$ by
the observation that
\begin{equation*}
M_w^{(1)}(s)=\int\limits_0^1\left(\log\Phi_w(x)-q_{r-1}(x,w)\right)x^{s-1}\,dx+
\sum_{\ell=k+1}^{r-1}\frac{b_\ell(w)}{s+\ell}.
\end{equation*}
The right hand side of this equation is a meromorphic function in $\Re s>-r$
and therefore provides the meromorphic continuation of $M_w^{(1)}(s)$. Thus
$M_w^{(1)}(s)$ has simple poles at $s=-\ell$ for $\ell\in\{k+1,k+2,\ldots\}$
with residues $b_\ell(w)$.

We proceed similarly for $M_w^{(2)}(s)$. This function is holomorphic for
$\Re s<-\rho$. From Theorem~\ref{thm1} we know that the periodic function $F$
has a uniformly convergent Fourier series
\begin{equation}\label{eq:Fourier}
F(u)=\sum_{m\in\mathbb{Z}}f_me^{2\pi imu}.
\end{equation}
From the fact that $F(u)$ is holomorphic for $|\Im u|<\beta\sigma$ (here and
later we set $\sigma=\frac1{\log\lambda}$) we actually get that
\begin{equation*}
f_m=\mathcal{O}_\eps(\exp(-2\pi\sigma(\beta-\eps)|m|))\quad\text{for all }
\eps>0.
\end{equation*}
Thus we can write
\begin{multline*}
M_w^{(2)}(s)=\int\limits_1^\infty\left(\log\Phi_w(x)-q_k(x,w)-
\sum_{m\in\mathbb{Z}}f_mx^{\rho+2\pi im\sigma}+
\frac{\log a_d}{d-1}+\log(-w)\right)x^{s-1}\,dx\\
+\sum_{m\in\mathbb{Z}}\frac{f_m}{s+\rho+2\pi im\sigma}-
\left(\frac{\log a_d}{d-1}+\log(-w)\right)\frac1{s},
\end{multline*}
where the right hand side is a meromorphic function in $\mathbb{C}$ with simple
poles at $s=-\rho-2\pi im\sigma$ ($m\in\mathbb{Z}$) with residue $f_m$ and at
$s=0$ with residue $-\frac{\log a_d}{d-1}-\log(-w)$. This provides the
meromorphic continuation of $M_w^{(2)}(s)$.

Thus we have proved
\begin{lemma}\label{lem1}
  The function $M_w(s)$ given by \eqref{eq:mellin} admits a meromorphic
  continuation to the whole complex plane with simple poles at $s=-\ell$ for
  $\ell\in\{k+1,k+2,\ldots\}$ (with $k=\fl{\rho}$), $s=0$, and $s=-\rho-2\pi
  im\sigma$ ($m\in\mathbb{Z}$) with $\sigma=\frac1{\log\lambda}$. The residues
  at these poles are given by (the constants $b_\ell(w)$ are given by
  \eqref{eq:logPhi})
\begin{equation}\label{eq:residues}
\begin{aligned}
\Res_{s=-\ell}M_w(s)&=b_\ell(w)\quad\text{for }\ell\in\{k+1,k+2,\ldots\}\\
\Res_{s=0}M_w(s)&=-\frac{\log a_d}{d-1}-\log(-w)\\
\Res_{s=-\rho-2\pi im\sigma}M_w(s)&=f_m\quad\text{for }m\in\mathbb{Z}.
\end{aligned}
\end{equation}
Under the conditions of Theorem~\ref{thm2}, Corollary~\ref{cor2}, or
Corollary~\ref{cor3} infinitely many values $f_m$ are non-zero.
\end{lemma}

Similar arguments applied to the function
\begin{equation}\label{eq:M0}
M_0(s)=\int_0^\infty
\left(\log\Phi(x)-\log x-\sum_{\ell=1}^k b_\ell x^\ell\right)x^{s-1}\,dx
\end{equation}
with $b_\ell$ given by
\begin{equation}\label{eq:logPhix}
\log\frac{\Phi(x)}x=\sum_{\ell=1}^\infty b_\ell x^\ell
\end{equation}
yield
\begin{lemma}\label{lem1b}
  The function $M_0(s)$ given by \eqref{eq:M0} admits a meromorphic
  continuation to the whole complex plane with simple poles at $s=-\ell$ for
  $\ell\in\{k+1,k+2,\ldots\}$ (with $k=\fl{\rho}$), and $s=-\rho-2\pi im\sigma$
  ($m\in\mathbb{Z}$) with $\sigma=\frac1{\log\lambda}$, and a double pole at
  $s=0$. The residues at these poles (the principal part resp.) are given by
  (the constants $b_\ell$ are given by \eqref{eq:logPhix})
\begin{equation}\label{eq:M0residues}
\begin{aligned}
\Res_{s=-\ell}M_0(s)&=b_\ell\quad\text{for }\ell\in\{k+1,k+2,\ldots\}\\
M_0(s)&=\frac1{s^2}-\frac{\log a_d}{d-1}\frac1s+\cdots\\
\Res_{s=-\rho-2\pi im\sigma}M_0(s)&=f_m\quad\text{for }m\in\mathbb{Z}.
\end{aligned}
\end{equation}
Under the conditions of Theorem~\ref{thm2}, Corollary~\ref{cor2}, or
Corollary~\ref{cor3} infinitely many values $f_m$ are non-zero.
\end{lemma}

In order to use the information derived for $M_w(s)$ for the meromorphic
continuation of the Dirichlet series
\begin{equation}\label{eq:zetaPhi}
\zeta_{\Phi,w}(s)=\sum_{\substack{\Phi(-\mu)=w\\\mu>0}}\mu^{-s}
\end{equation}
we derive an alternative expression for $M_w(s)$ based on \eqref{eq:hadamard}
\begin{equation}\label{eq:tildeM}
\begin{aligned}
M_w(s)&=\int\limits_0^\infty\left(\log\Phi_w(x)-q_k(x,w)\right)x^{s-1}\,dx\\
&=\sum_{\Phi_w(-\mu)=0}\int\limits_0^\infty
\left(\log\left(1+\frac x\mu\right)-\frac x\mu+\frac{x^2}{2\mu^2}+\cdots
+(-1)^{k}\frac{x^k}{k\mu^k}\right)x^{s-1}\,dx.
\end{aligned}
\end{equation}
Furthermore, the second integral in \eqref{eq:tildeM} can
be evaluated as
\begin{equation}\label{eq:zetaPhi-M}
M_w(s)=\zeta_{\Phi,w}(-s)\frac\pi{s\sin\pi s}
\end{equation}
by the fact that
\begin{equation*}
\int\limits_0^\infty\left(\log(1+x)-x+\frac{x^2}2+\cdots+
(-1)^k\frac{x^k}k\right)x^{s-1}\,dx=\frac\pi{s\sin\pi s}\quad\text{for }
-k-1<\Re s<-k
\end{equation*}
and the fact that the right hand side is the Mellin transform of a ``harmonic
sum'' (cf.~\cite{Oberhettinger1974:tables_mellin_transforms,
Paris_Kaminski2001:asymptotics_mellin_barnes}).
Equation \eqref{eq:zetaPhi-M} provides us with the meromorphic continuation of
$\zeta_{\Phi,w}(s)$. This is somewhat a reversion of the ideas used in
\cite{Mellin1903:dirichletsche_reihen}.

Summing up, we have obtained
\begin{theorem}\label{thm3}
  Let $\Phi$ be the solution of the functional equation \eqref{eq:schr}. Then
  the Dirichlet generating function $\zeta_{\Phi,w}(s)$ of the solutions of
  $\Phi(-\mu)=w$ for $w<0$ given in \eqref{eq:zetaPhi} admits a meromorphic
  continuation to the whole complex plane. Under the conditions of
  Theorem~\ref{thm2}, Corollary~\ref{cor2}, or Corollary~\ref{cor3} there exist
  infinitely many simple poles at points of the form $s=\rho+2\pi im\sigma$
  ($\sigma=\frac1{\log\lambda}$, $m\in\Z$).  The following special values are
  known ($\rho=\frac{\log d}{\log\lambda}$ and $f_m$ are given by
  \eqref{eq:Fourier})
\begin{align*}
\Res_{\rho+2\pi im\sigma}\zeta_{\Phi,w}(s)&=
-\frac{f_m}\pi\left(\rho+2\pi im\sigma\right)
\sin\pi\left(\rho+2\pi im\sigma\right)\text{ for }m\in\Z\\
\zeta_{\Phi,w}(m)&=0\quad\text{for }m<\rho, m\in\mathbb{Z}\\
\zeta_{\Phi,w}'(0)&=\frac{\log a_d}{d-1}+\log(-w)\\
\zeta_{\Phi,w}(m)&=(-1)^{m-1}mb_m(w)\quad\text{for }m\in\mathbb{N},
\end{align*}
where $b_m(w)$ is given by \eqref{eq:logPhi}.
\end{theorem}

Similar arguments applied to $\Phi(x)/x$ yield
\begin{theorem}\label{thm4}
  Let $\Phi$ be the solution of the functional equation \eqref{eq:schr}. Then
  the Dirichlet generating function $\zeta_{\Phi,0}(s)$ of the solutions of
  $\Phi(-\mu)=0$ for $\mu>0$ given in \eqref{eq:zetaPhi} admits a meromorphic
  continuation to the whole complex plane.  Under the conditions of
  Theorem~\ref{thm2}, Corollary~\ref{cor2}, or Corollary~\ref{cor3} there exist
  infinitely many simple poles at points of the form $s=\rho+2\pi im\sigma$
  ($\sigma=\frac1{\log\lambda}$, $m\in\Z$).  The following special values are
  known ($\rho=\frac{\log d}{\log\lambda}$)
\begin{align*}
\Res_{\rho+2\pi im\sigma}\zeta_{\Phi,0}(s)&=
-\frac{f_m}\pi\left(\rho+2\pi im\sigma\right)
\sin\pi\left(\rho+2\pi im\sigma\right)\text{ for }m\in\Z\\
\zeta_{\Phi,0}(m)&=0\quad\text{for }m<\rho, m\in\mathbb{Z}\setminus\{0\}\\
\zeta_{\Phi,0}(0)&=1\\
\zeta_{\Phi,0}'(0)&=\frac{\log a_d}{d-1}\\
\zeta_{\Phi,0}(m)&=(-1)^{m-1}mb_m\quad\text{for }m\in\mathbb{N},
\end{align*}
where $b_m$ is given by \eqref{eq:logPhix}.
\end{theorem}

\begin{remark}\label{rem5}
By Lemma~\ref{lem5} we have
\begin{equation*}
\zeta_{\Phi,w}(s)=\lim_{n\to\infty}
\sum_{z\in p^{(-n)}(\{w\})}(\lambda^nz)^{-s},
\end{equation*}
which is the function $\zeta_{p,w}(2s)$ studied in
\cite{Teplyaev2004:spectral_zeta_function,
Teplyaev2005:spectral_zeta_functions}.
\end{remark}

\section{Meromorphic continuation of the zeta-function}
\label{sec:analyt-cont-zeta}
Up to now we have only considered the eigenvalues originating from $m=0$ in
\eqref{eq:decimation} for fixed $w\in A$. In general the multiplicity of
eigenvalues for $m>0$ is given by a linear recurrent sequence $\beta_m(w)$.
This sequence has a rational generating function
\begin{equation*}
B_w(x)=\sum_{m=0}^\infty \beta_m(w) x^m=\frac{P_w(x)}{Q_w(x)}
\end{equation*}
with $P_w,Q_w\in\Z[x]$. Let $r_w$ denote the radius of convergence of $B_w$,
then by Pringsheim's theorem $x=r_w$ is a pole of $B_w$; let $k_w$ denote the
order of this pole. Then the positivity of $\beta_m(w)$ implies that the
coefficient $c_{k_w}(w)$ in the Laurent-expansion around $x=r_w$
\begin{equation*}
B_w(x)=\frac{c_{k_w}(w)}{\left(1-\frac x{r_w}\right)^{k_w}}+\cdots
\end{equation*}
is positive.

Putting everything together, we can write
\begin{equation}\label{eq:zeta-final}
\zeta_\Delta(s)=\sum_{-\Delta u=\mu u}\frac1{\mu^s}=
\sum_{w\in A}\frac{P_w(\lambda^{-s})}{Q_w(\lambda^{-s})}\zeta_{\Phi,w}(s),
\end{equation}
which provides the meromorphic continuation of $\zeta_\Delta$ to the whole
complex plane. If $\log_\lambda d<\frac12d_S$ then all the functions
$\zeta_{\Phi,w}(s)$ are holomorphic in $\Re s>\frac12d_S-\eps$ for some
$\eps>0$. Since $\zeta_\Delta(s)$ has a simple pole at $s=\frac12d_S$ by the
fact that $L(x)\asymp x^{\frac12d_S}$, at least
one of the rational functions $B_w(x)$ has to have a pole at
$x=\lambda^{-\frac12d_S}$. All poles of functions $B_w$ at
$x=\lambda^{-\frac12d_S}$ have to be simple. Let $W$ denote the set of those
$w$ with $r_w=\lambda^{-\frac12d_S}$.  By our observation on the sign of
$c_1(w)$ for $w\in W$ the Dirichlet series
\begin{equation*}
\eta(s)=-\sum_{w\in W}c_1(w)\zeta_{\Phi,w}(s)
\end{equation*}
has positive coefficients.  We apply
\cite[Theorem~9.5, p.~184]{Lapidus_Frankenhuysen2000:fractal_geometry_number}
to the function $\eta(s)$ to see that $\eta(\frac12d_S+ik\tau)=0$ cannot hold
for fixed $\tau$ and all $k\in\Z\setminus\{0\}$. Thus the function
\begin{equation*}
\sum_{w\in W}
B_w(\lambda^{-s})\zeta_{\Phi,w}(s)
\end{equation*}
has a simple pole at $s=\frac12d_S$ and at least two non-real poles on the line
$\Re s=\frac12d_S$. Since the remaining summands in \eqref{eq:zeta-final} have
no poles on this line, these points remain poles of $\zeta_\Delta(s)$.

Summing up, we have proved:
\begin{theorem}\label{thm8}
  Let $G$ be a p.~c.~f.\ self-similar compact fractal, whose Laplace operator
  $\Delta$ admits spectral decimation in the sense of Definition~\ref{def1}
  with a polynomial of degree $d$. Then the Dirichlet generating function of
  the eigenvalues of $\Delta$
\begin{equation*}
\zeta_\Delta(s)=\sum_{-\Delta u=\mu u}\frac1{\mu^s},
\end{equation*}
has a meromorphic continuation to the whole complex plane with poles contained
in a finite union of sets $\{\rho_k+2\pi im\sigma\mid m\in\mathbb{Z}\}$, where
$\sigma=\frac1{\log\lambda}$ and $\lambda$ is the parameter coming from
spectral decimation. There is a simple pole at $s=\frac12d_S$ ($d_S$ denoting
the spectral dimension of $G$). If $\log_\lambda d<\frac12d_S$ then
$\zeta_\Delta(s)$ has at least two non-real poles on the line
$\Re s=\frac12d_S$.
\end{theorem}
\begin{remark}\label{rem6}
The case of $G=[0,1]$ which gives the Riemann zeta function and has
$\log_\lambda d=\frac12d_S$ shows that the condition $\log_\lambda
d<\frac12d_S$ is needed for the last assertion.
\end{remark}
\begin{theorem}\label{thm9}
  Let $G$ be a p.~c.~f.\ self-similar compact fractal, whose Laplace operator
  $\Delta$ admits spectral decimation in the sense of Definition~\ref{def1}.
  Then the following are equivalent:
  \begin{enumerate}
  \item\label{enum1} $\zeta_\Delta(s)$ has at least two non-real poles in the
    set $\frac12d_S+2\pi i\sigma\Z$ ($\sigma=\frac1{\log\lambda}$).
  \item\label{enum2} the limit $\lim_{x\to\infty}x^{-\frac12 d_S}L(x)$ does not
    exist, where $L(x)$ denotes the eigenvalue counting function
    \eqref{eq:count}
  \item\label{enum3} the limit $\lim_{t\to0+}P(t)t^{\frac12d_S}$ does not
    exist, where $P(t)$ denotes the trace of the heat kernel \eqref{eq:heat}.
  \end{enumerate}
\end{theorem}
\begin{proof}
  ``\eqref{enum1}$\Leftrightarrow$\eqref{enum2}'': the general theory of
  Dirichlet series (cf.~\cite{Hardy_Riesz1964:general_theory_dirichlet})
  implies that $\zeta_\Delta(s)(\lambda^s-\lambda^{\frac12d_S})$ is of
  polynomial growth along vertical lines. Thus for some $k\in\N$ the last
  integral in the Mellin-Perron formula
  (cf.~\cite{Tenenbaum1995:introduction_to_analytic})
\begin{multline*}
\sum_{\mu<x}\left(1-\frac\mu x\right)^k=
\frac1{2\pi i}\int\limits_{\frac12d_S+1-i\infty}^{\frac12d_S+1+i\infty}
\zeta_\Delta(s)x^s\frac{ds}{s(s+1)\cdots(s+k)}=\\
\sum_{\ell\in\Z}\Res_{s=\frac12d_S+2\pi i\ell\sigma}
\zeta_\Delta(s)x^s\frac1{s(s+1)\cdots(s+k)}+
\frac1{2\pi i}\int\limits_{\frac12d_S-\eps-i\infty}^{\frac12d_S-\eps+i\infty}
\zeta_\Delta(s)x^s\frac{ds}{s(s+1)\cdots(s+k)}
\end{multline*}
converges for some $\eps>0$. The last integral is
$\BigOh(x^{\frac12d_S-\eps})$. The series constitutes $x^{\frac12d_S}$ times a
Fourier series in $\log_\lambda x$, which is constant, if and only if there are
no non-real poles on the line $\Re s=\frac12d_S$. Since the left-hand side of
the equation is an iterated integral of $L(x)$, $x^{-\frac12d_S}L(x)$ shows
oscillations, if and only if the right-hand side does.

 ``\eqref{enum1}$\Leftrightarrow$\eqref{enum3}'': again by Mellin inversion we
 get
\begin{multline*}
P(t)=\frac1{2\pi i}\int\limits_{\frac12d_S+1-i\infty}^{\frac12d_S+1+i\infty}
\zeta_\Delta(s)\Gamma(s)t^{-s}\,ds=\\
\sum_{\ell\in\Z}\Res_{s=\frac12d_S+2\pi i\ell\sigma}\zeta_\Delta(s)\Gamma(s)t^{-s}+
\frac1{2\pi i}\int\limits_{\frac12d_S-\eps-i\infty}^{\frac12d_S-\eps+i\infty}
\zeta_\Delta(s)\Gamma(s)t^{-s}\,ds.
\end{multline*}
Again the series is $t^{-\frac12d_S}$ times a Fourier series in $\log_\lambda
t$ and the last integral constitutes an error term
$\BigOh(t^{-\frac12d_S+\eps})$. The same argument as above proves the presence
of oscillations, if and only if $\zeta_\Delta$ has non-real poles on the line
$\Re s=\frac12d_S$.
\end{proof}
\begin{remark}\label{rem4}
  In \cite[p.~105]{Kigami_Lapidus1993:weyl_problem_spectral} it is conjectured
  that the limit $\lim_{x\to\infty}x^{-\frac12 d_S}L(x)$ does not exist for
  non-integer $d_S$ in the ``\emph{lattice case}''. The non-existence of the
  limit has been confirmed in the case of existence of localised eigenfunctions
  (cf.~\cite{Barlow_Kigami1997:localized_eigenfunctions_laplacian,
    Malozemov_Teplyaev1995:pure_point_spectrum, Sabot2000:pure_point_spectrum})
  and in the case of the existence of spectral gaps
  (cf.~\cite{Strichartz2005:laplacians_on_fractals}).
\end{remark}
\begin{remark}\label{rem7}
  In \cite{Kroen_Teufl2004:asymptotics_transition_probabilities} it has been
  proved that the $n$-step transition probabilities $p_n(x,y)$ (the discrete
  analogue of the heat-kernel) on certain self-similar graphs satisfies an
  asymptotic relation
\begin{equation*}
p_n(x,y)\sim n^{-\frac12d_S}F(\log_\lambda n)
\end{equation*}
with a periodic function $F$. It is shown that $F$ is not constant, if the
Julia set of a related rational function is a Cantor set. In
\cite{Teufl2006:asymptotic_behaviour_analytic} it is proved for the same class
of graphs that $F$ is not constant, if $d_S$ is irrational.
\end{remark}
\begin{remark}\label{rem8}
  Theorems~\ref{thm8} and~\ref{thm9} together confirm the conjecture
  (cf.~\cite[p.~105]{Kigami_Lapidus1993:weyl_problem_spectral}) that the limit
  $\lim_{x\to\infty}x^{-\frac12 d_S}L(x)$ does not exist for fractals admitting
  spectral decimation in the sense of Definition~\ref{def1} with a polynomial
  of degree $d$ and $\log_\lambda d<\frac12d_S$.
\end{remark}
\section{Examples}\label{sec:examples}
\begin{example}\label{ex1}
  In \cite{Fukushima_Shima1992:spectral_analysis_sierpinski,
    Shima1991:eigenvalue_problems_random} the eigenvalues of the self-similar
  Dirichlet-Laplacian on the $K$-dimensional Sierpi\'nski gasket have been
  identified. The $K$-dimensional gasket admits spectral decimation with the
  polynomial $p_K(x)=x(K+3+x)$. The eigenvalues of $\Delta$ are derived from
  the preiterates of the elements of the set $A=\{-2,-(K+1),-(K+3)\}$ with
  multiplicities
\begin{align*}
  \beta_m(-2)&=
  \begin{cases}
  1&\text{for }m=1\\
  0&\text{otherwise,}
  \end{cases}\\
  \beta_m(-(K+1))&=\frac{K-1}2(K+1)^{m-1}-\frac{K+1}2\text{ for }m\geq 2,\\
  \beta_m(-(K+3))&=\frac{K-1}2(K+1)^{m-1}+\frac{K+1}2\text{ for }m\geq 1.
\end{align*}
This yields (with $\lambda=K+3$)
\begin{multline}\label{eq:ex1}
\zeta_\Delta(s)=\lambda^{-s}\zeta_{\Phi,-2}(s)+
\zeta_{\Phi,-(K+1)}(s)\sum_{m=2}^\infty \beta_{m}(-(K+1))\lambda^{-ms}+\\
\zeta_{\Phi,-(K+3)}(s)\sum_{m=1}^\infty\beta_m(-(K+3))\lambda^{-ms}=\\
\lambda^{-s}\zeta_{\Phi,-2}(s)+
\left(\frac{K^2-1}{2(\lambda^s-(K+1))}-
\frac{K+1}{2(\lambda^s-1)}\right)\lambda^{-s}\zeta_{\Phi,-(K+1)}(s)+\\
\left(\frac{K-1}{2(\lambda^s-(K+1))}+
\frac{K+1}{2(\lambda^s-1)}\right)\zeta_{\Phi,-(K+3)}(s).
\end{multline}
We first observe that the poles of $\zeta_{\Phi,\cdot}(s)$ at $s=\rho+2\pi
im\sigma$, which are given by $\lambda^s=2$, cancel by the fact that the
residues of $\zeta_{\Phi,w}$ at these points are independent of $w$ and the
values of the rational factors sum up to $0$. Furthermore, we observe
\begin{equation*}
\zeta_{\Phi,-(K+3)}(s)=\left(\lambda^s-1\right)\zeta_{\Phi,0}(s)
\end{equation*}
by the fact that $\Phi(z)=-(K+3)$, if and only if $\Phi(\lambda z)=0$ and
$\Phi(z)\neq0$. This implies that the poles at $s=2\pi im\sigma$ of the
rational factor of the last summand in \eqref{eq:ex1} are cancelled.
Using Mellin transform, we see that the cancellation of poles on the line
$\Re s=0$ is equivalent to the existence of the limit
\begin{equation*}
\lim_{x\to\infty}\frac{x^{\log_\lambda(-w)}}{\Phi(x)}
\prod_{n=1}^\infty\left(1-\frac1w\Phi(x)\right)
\end{equation*}
We did some numerical computations, which indicate that there are oscillations
for $w=-3$ and $K=2$.  Nevertheless, the pole of the rational factor at $s=0$
is cancelled by the zero of $\zeta_{\Phi,-K-1}$.

From this we
derive, for instance,
\begin{align*}
\zeta_\Delta(0)&=\frac{K+1}2\left(1-2\frac{\log(K+1)}{\log(K+3)}\right)\\
\zeta'_\Delta(0)&=\log2+\frac{2K^2+3K+1}{2K}\log(K+1)-
\frac{K^2+3K-2}{4K}\log(K+3)\\
&+\frac{K+1}{4\log(K+3)}\zeta''_{\Phi,-K-3}(0)-
\frac{K+1}{2\log(K+3)}\zeta''_{\Phi,-K-1}(0)\\
\zeta_\Delta(1)&=\frac{K^2+3K-1}{2(K+2)(K+3)}
\end{align*}
\end{example}
\begin{example}\label{ex2}
  In \cite{Teplyaev1998:spectral_analysis_gaskets,
    Teplyaev2004:spectral_zeta_function,Teplyaev2005:spectral_zeta_functions}
  the spectrum of the Neumann-Laplacian on the $2$-dimensional Sierpi\'nski
  gasket has been studied. It has been shown that it admits spectral decimation
  with the polynomial $p(x)=x(5+x)$. The Neumann-eigenvalues are derived from
  the preiterates of the set $A=\{-3,-5\}$ with multiplicities
\begin{align*}
  \beta_m(-3)&=\frac{3^{m-1}+3}2\text{ for }m\geq1\\
  \beta_m(-5)&=\frac{3^{m-1}-1}2\text{ for }m\geq2.
\end{align*}
This gives
\begin{multline*}
\zeta_{\Delta}(s)=\zeta_{\Phi,-3}(s)\sum_{m=1}^\infty\beta_m(-3)5^{-ms}+
\zeta_{\Phi,-5}(s)\sum_{m=2}^\infty\beta_m(-5)5^{-ms}=\\
\left(\frac1{2(5^s-3)}+\frac3{2(5^s-1)}\right)\zeta_{\Phi,-3}(s)+
\left(\frac3{2(5^s-3)}-\frac1{2(5^s-1)}\right)5^{-s}\zeta_{\Phi,-5}(s),
\end{multline*}
which is in accordance with
\cite{Teplyaev2004:spectral_zeta_function,
Teplyaev2005:spectral_zeta_functions}. Notice that there is no pole at $s=0$,
the poles at the solutions of $5^s=1$ are cancelled in the second summand by
the observation in Example~\ref{ex1}. Numerical experiments seem to indicate
that the poles of the first summand on the imaginary axis do not cancel with
zeros of $\zeta_{\Phi,-3}$.

By the same arguments as in Example~\ref{ex1} the poles at the solutions
of $5^s=2$ are cancelled. Furthermore, in the second summand the poles at the
solutions of $5^s=1$ are cancelled by zeros of $\zeta_{\Phi,-5}$. The poles of
the rational factors at $s=0$ cancel and we obtain
\begin{align*}
  \zeta_{\Delta}(0)&=\frac32\log_53-\frac12\\
\zeta'_{\Delta}(0)&=\frac{3\zeta''_{\Phi,-3}(0)-\zeta''_{\Phi,-5}(0)}{4\log5}-
\log3\\
\zeta_{\Delta}(1)&=\frac7{30}\\
\zeta_{\Delta}(2)&=\frac1{150}
\end{align*}

The second derivatives $\zeta''_{\Phi,w}(0)$ can be computed numerically by the
following observations. Around $s=0$ we can write
$M_w(s)=-\frac{\log(-w)}s+H_w(s)$ with a function $H_w$ holomorphic around $0$,
and by \eqref{eq:zetaPhi-M} we have $2H_w(0)=\zeta''_{\Phi,-w}(0)$. We have
\begin{multline*}
(1-2\cdot 5^s)H_w(s)+d\frac{\log(-w)}s(5^s-1)=\\
\int_0^1\log\left(\frac{1-\frac1w\Phi(x)}{(1-\frac1w\Phi(\frac x5))^2}\right)
x^{s-1}\,dx+
\int_1^\infty
\log\left(\frac{1-\frac1w\Phi(x)}{-w(1-\frac1w\Phi(\frac x5))^2}\right)
x^{s-1}\,dx,
\end{multline*}
and setting $s=0$
\begin{multline*}
H_w(0)=2\log(-w)\log5-
\int\limits_0^1\!\!
\log\left(\frac{1-\frac1w\Phi(x)}{(1-\frac1w\Phi(\frac x5))^2}\right)
\,\frac{dx}x-\int\limits_1^\infty\!\!
\log\left(\frac{1-\frac1w\Phi(x)}{-w(1-\frac1w\Phi(\frac x5))^2}\right)
\,\frac{dx}x.
\end{multline*}
The first integral can be computed as the rapidly convergent series
\begin{equation*}
\int_0^1\log\left(\frac{1-\frac1w\Phi(x)}{(1-\frac1w\Phi(\frac x5))^2}\right)
\,\frac{dx}x=\sum_{n=1}^\infty\frac{b_n(w)}n(1-2\cdot5^{-n}).
\end{equation*}
The power series for $\log(1-\frac1w\Phi(z))$ has radius of convergence the
modulus of the smallest solution of $\Phi(z)=w$, which is much larger than $1$
in the cases of interest. This gives exponential convergence.

For computing the second interval we observe that $\Phi(x)\geq\exp(0.9 x^\rho)$
for $x\geq5$, which can be shown by discussing $\Phi$ on the finite interval
$5\leq x\leq25$ and the extending by the obvious inequality
$\Phi(5x)\geq\Phi(x)^2$. This together with an inequality for the logarithm
yields (for $T\geq25$)
\begin{equation*}
0\geq\int_{T}^\infty\log
\left(\frac{1-\frac1w\Phi(x)}{-w(1-\frac1w\Phi(\frac x5))^2}\right)
\,\frac{dx}x\geq(2w-1)\int_T^\infty\exp(-0.9x^\rho)\frac{dx}x.
\end{equation*}
Thus the improper integral can be computed by truncation at a finite value and
estimating the remainder integral as above.

Using these ideas we computed (using \texttt{Mathematica})
\begin{align*}
H_{-3}(0)&=5.23995\ 51500\ldots\\
H_{-5}(0)&=9.06601\ 63789\ldots\\
\zeta'_\Delta(0)&=0.96852\ 21499\ldots.
\end{align*}
If there were no poles of $\zeta_\Delta$ on the line $\Re s=0$,
$\exp(-\zeta'_\Delta(0))$ were the value of the regularised product of the
eigenvalues (cf.\cite{Jorgenson_Lang1993:basic_analysis_regularized}), or the
Fredholm determinant of $\Delta$.
\end{example}
\begin{ackno}
The second author is indebted to Michel Lapidus and Alexander Teplyaev for
helpful discussions. The authors are grateful to an anonymous referee for
valuable suggestions.
\end{ackno}

\providecommand{\bysame}{\leavevmode\hbox to3em{\hrulefill}\thinspace}
\providecommand{\MR}{\relax\ifhmode\unskip\space\fi MR }
\providecommand{\MRhref}[2]{%
  \href{http://www.ams.org/mathscinet-getitem?mr=#1}{#2}
}
\providecommand{\href}[2]{#2}


\begin{thebibliography}{10}

\bibitem{Barlow1998:diffusions_on_fractals}
M.~T. Barlow, \emph{Diffusions on fractals}, Lectures on probability theory and
  statistics (Saint-Flour, 1995), Springer Verlag, Berlin, 1998, pp.~1--121.

\bibitem{Barlow_Kigami1997:localized_eigenfunctions_laplacian}
M.~T. Barlow and J.~Kigami, \emph{Localized eigenfunctions of the {L}aplacian
  on p.c.f.\ self-similar sets}, J. London Math. Soc. (2) \textbf{56} (1997),
  320--332.

\bibitem{Barlow_Perkins1988:brownian_motion_sierpinski}
M.~T. Barlow and E.~A. Perkins, \emph{Brownian motion on the {S}ierpi\'nski
  gasket}, Probab. Theory Relat. Fields \textbf{79} (1988), 543--623.

\bibitem{Beardon1991:iteration_rational_functions}
A.~F. Beardon, \emph{Iteration of rational functions}, Graduate Texts in
  Mathematics, no. 132, Springer Verlag, 1991.

\bibitem{Biggins_Bingham1991:near_constancy_phenomena}
J.~D. Biggins and N.~H. Bingham, \emph{Near-constancy phenomena in branching
  processes}, Math. Proc. Camb. Philos. Soc. \textbf{110} (1991), 545--558.

\bibitem{Boas1954:entire_functions}
R.~P. Boas, Jr., \emph{Entire functions}, Academic Press Inc., New York, 1954.

\bibitem{Doetsch1971:handbuch_der_laplace}
G.~Doetsch, \emph{Handbuch der {L}aplace-{T}ransformation. {B}and {I}:
  {T}heorie der {L}aplace-{T}ransformation}, Birkh\"auser Verlag, Basel, 1971,
  Verbesserter Nachdruck der ersten Auflage 1950, Lehrb\"ucher und Monographien
  aus dem Gebiete der exakten Wissenschaften. Mathematische Reihe, Band 14.

\bibitem{Dubuc1982:etude_theorique_numerique}
S.~Dubuc, \emph{Etude th\'eorique et num\'erique de la fonction de
  {K}arlin-{M}c{G}regor}, J. Analyse Math. \textbf{42} (1982), 15--37.

\bibitem{Fukushima_Shima1992:spectral_analysis_sierpinski}
M.~Fukushima and T.~Shima, \emph{On a spectral analysis for the {S}ierpi\'nski
  gasket}, Potential Anal. \textbf{1} (1992), 1--35.

\bibitem{Grabner1997:functional_iterations_stopping}
P.~J. Grabner, \emph{Functional iterations and stopping times for {B}rownian
  motion on the {S}ierpi\'nski gasket}, Mathematika \textbf{44} (1997),
  374--400.

\bibitem{Hardy_Riesz1964:general_theory_dirichlet}
G.~H. Hardy and M.~Riesz, \emph{The general theory of {D}irichlet's series},
  Cambridge Tracts in Mathematics and Mathematical Physics, No. 18,
  Stechert-Hafner, Inc., New York, 1964.

\bibitem{Jorgenson_Lang1993:basic_analysis_regularized}
J.~Jorgenson and S.~Lang, \emph{Basic analysis of regularized series and
  products}, Lecture Notes in Mathematics, vol. 1564, Springer-Verlag, Berlin,
  1993.

\bibitem{Karlin_Mcgregor1968:embeddability_discrete_time}
S.~Karlin and J.~McGregor, \emph{Embeddability of discrete time simple
  branching processes into continuous time branching processes}, Trans. Amer.
  Math. Soc. \textbf{132} (1968), 115--136.

\bibitem{Kigami1993:harmonic_calculus_p}
J.~Kigami, \emph{A harmonic calculus for p.c.f. self-similar sets}, Trans.
  Amer. Math. Soc. \textbf{335} (1993), 721--755.

\bibitem{Kigami1998:distributions_localized_eigenvalues}
J.~Kigami, \emph{Distributions of localized eigenvalues of {L}aplacians on post
  critically finite self-similar sets}, J. Funct. Anal. \textbf{156} (1998),
  170--198.

\bibitem{Kigami2001:analysis_fractals}
\bysame, \emph{Analysis on fractals}, Cambridge Tracts in Mathematics, vol.
  143, Cambridge University Press, Cambridge, 2001.

\bibitem{Kigami_Lapidus1993:weyl_problem_spectral}
J.~Kigami and M.~L. Lapidus, \emph{Weyl's problem for the spectral distribution
  of {L}aplacians on p.c.f.\ self-similar fractals}, Comm. Math. Phys.
  \textbf{158} (1993), 93--125.

\bibitem{Kroen2002:green_functions_self}
B.~Kr{\"o}n, \emph{Green functions on self-similar graphs and bounds for the
  spectrum of the {L}aplacian}, Ann. Inst. Fourier (Grenoble) \textbf{52}
  (2002), 1875--1900.

\bibitem{Kroen_Teufl2004:asymptotics_transition_probabilities}
B.~Kr{\"o}n and E.~Teufl, \emph{Asymptotics of the transition probabilities of
  the simple random walk on self-similar graphs}, Trans. Amer. Math. Soc.
  \textbf{356} (2004), 393--414.

\bibitem{Kuczma1963:schroeder_equation}
M.~Kuczma, \emph{On the {S}chr\"oder equation}, Rozprawy Mat. \textbf{34}
  (1963), 50.

\bibitem{Lapidus1994:analysis_fractals_laplacians}
M.~L. Lapidus, \emph{Analysis on fractals, {L}aplacians on self-similar sets,
  noncommutative geometry and spectral dimensions}, Topol. Methods Nonlinear
  Anal. \textbf{4} (1994), 137--195.

\bibitem{Lapidus_Frankenhuysen2000:fractal_geometry_number}
M.~L. Lapidus and M.~van Frankenhuysen, \emph{Fractal geometry and number
  theory}, Birkh\"auser Boston Inc., Boston, MA, 2000, Complex dimensions of
  fractal strings and zeros of zeta functions.

\bibitem{Lindstroem1990:brownian_motion_nested}
T.~Lindstr{\o}m, \emph{Brownian motion on nested fractals}, Mem. Amer. Math.
  Soc., vol. 420, Amer. Math. Soc., 1990.

\bibitem{Malozemov_Teplyaev1995:pure_point_spectrum}
L.~Malozemov and A.~Teplyaev, \emph{Pure point spectrum of the {L}aplacians on
  fractal graphs}, J. Funct. Anal. \textbf{129} (1995), 390--405.

\bibitem{Malozemov_Teplyaev2003:self_similarity_operators}
\bysame, \emph{Self-similarity, operators and dynamics}, Math. Phys. Anal.
  Geom. \textbf{6} (2003), 201--218.

\bibitem{Mellin1903:dirichletsche_reihen}
H.~Mellin, \emph{Die {D}irichlet'schen {R}eihen, die zahlentheoretischen
  {F}unktionen und die unendlichen {P}rodukte von endlichem {G}eschlecht}, Acta
  Math. \textbf{28} (1903), 37--64.

\bibitem{Minakshisundaram_Pleijel1949:some_properties_eigenfunctions}
S.~Minakshisundaram and {\AA}.~Pleijel, \emph{Some properties of the
  eigenfunctions of the {L}aplace-operator on {R}iemannian manifolds}, Canadian
  J. Math. \textbf{1} (1949), 242--256.

\bibitem{Oberhettinger1974:tables_mellin_transforms}
F.~Oberhettinger, \emph{Tables of {M}ellin transforms}, Springer-Verlag, New
  York, 1974.

\bibitem{Paris_Kaminski2001:asymptotics_mellin_barnes}
R.~B. Paris and D.~Kaminski, \emph{Asymptotics and {M}ellin-{B}arnes
  integrals}, Encyclopedia of Mathematics and its Applications, vol.~85,
  Cambridge University Press, Cambridge, 2001.

\bibitem{Rammal_Toulouse1983:random_walks_fractal}
R.~Rammal and G.~Toulouse, \emph{Random walks on fractal structures and
  percolation clusters}, J. Physique Lettres \textbf{44} (1983), L13--L22.

\bibitem{Rosenberg1997:laplacian_riemannian_manifold}
S.~Rosenberg, \emph{The {L}aplacian on a {R}iemannian manifold}, London
  Mathematical Society Student Texts, vol.~31, Cambridge University Press,
  Cambridge, 1997, An introduction to analysis on manifolds.

\bibitem{Sabot2000:pure_point_spectrum}
C.~Sabot, \emph{Pure point spectrum for the {L}aplacian on unbounded nested
  fractals}, J. Funct. Anal. \textbf{173} (2000), 497--524.

\bibitem{Shima1991:eigenvalue_problems_random}
T.~Shima, \emph{On eigenvalue problems for the random walk on the
  {S}ierpi\'nski pre-gaskets}, Japan J. Appl. Ind. Math. \textbf{8} (1991),
  127--141.

\bibitem{Shima1993:eigenvalue_problem_laplacian}
T.~Shima, \emph{The eigenvalue problem for the {L}aplacian on the
  {S}ierpi\'nski gasket}, Asymptotic problems in probability theory: stochastic
  models and diffusions on fractals (Sanda/Kyoto, 1990), Pitman Res. Notes
  Math. Ser., vol. 283, Longman Sci. Tech., Harlow, 1993, pp.~279--288.

\bibitem{Shima1996:eigenvalue_problems_laplacians}
\bysame, \emph{On eigenvalue problems for {L}aplacians on p.c.f. self-similar
  sets}, Japan J. Indust. Appl. Math. \textbf{13} (1996), 1--23.

\bibitem{Strichartz1999:some_properties_laplacians}
R.~S. Strichartz, \emph{Some properties of {L}aplacians on fractals}, J. Funct.
  Anal. \textbf{164} (1999), 181--208.

\bibitem{Strichartz2003:fractafolds_spectra}
\bysame, \emph{Fractafolds based on the {S}ierpi\'nski gasket and their
  spectra}, Trans. Amer. Math. Soc. \textbf{355} (2003), 4019--4043.

\bibitem{Strichartz2005:laplacians_on_fractals}
\bysame, \emph{Laplacians on fractals with spectral gaps have nicer {F}ourier
  series}, Math. Res. Lett. \textbf{12} (2005), 269--274.

\bibitem{Tenenbaum1995:introduction_to_analytic}
G.~Tenenbaum, \emph{Introduction to analytic and probabilistic number theory},
  Cambridge Studies in Advanced Mathematics, no.~46, Cambridge University
  Press, 1995.

\bibitem{Teplyaev1998:spectral_analysis_gaskets}
A.~Teplyaev, \emph{Spectral analysis on infinite {S}ierpi\'nski gaskets}, J.
  Funct. Anal. \textbf{159} (1998), 537--567.

\bibitem{Teplyaev2004:spectral_zeta_function}
\bysame, \emph{Spectral zeta function of symmetric fractals}, Fractal geometry
  and stochastics III (C.~Bandt, U.~Mosco, and M.~Z{\"a}hle, eds.), Progr.
  Probab., vol.~57, Birkh\"auser, Basel, 2004, pp.~245--262.

\bibitem{Teplyaev2005:spectral_zeta_functions}
\bysame, \emph{Spectral zeta functions of fractals and the complex dynamics of
  polynomials}, available at \url{http://arxiv.org/pdf/math.SP/0505546}, 2005.

\bibitem{Teufl2006:asymptotic_behaviour_analytic}
E.~Teufl, \emph{On the asymptotic behaviour of analytic solutions of linear
  iterative functional equations}, Aequationes Math. (2006), to appear.

\bibitem{Valiron1954:fonctions_analytiques}
G.~Valiron, \emph{Fonctions analytiques}, Presses Universitaires de France,
  Paris, 1954.

\end{thebibliography}
\end{document}